\def\ifl{independence friendly logic}

\def\proof{\noindent{\bf Proof.}\hspace{3mm}}
\def\qed{\hfill{\bf Q.E.D.}\medskip}

\newcommand{\rk}{{\rm w}}

\def\={\mathop{=}}

\def\p{\phi}
\def\s{\psi}
\def\h{\theta}

\def\fl{first order logic}
\newtheorem{theorem}{Theorem}
\newtheorem{definition}[theorem]{Definition}
\newtheorem{lemma}[theorem]{Lemma}

\newtheorem{corollary}[theorem]{Corollary}

\def\dom{\mathop{\rm dom}}
\def\ran{\mathop{\rm ran}}

\def\R{{\Bbb R}}
\newcommand{\open}{\Bbb}

\newcommand{\oN}{{\open N}}

\def\ef{{\rm EF}}

\def\EF{Ehrenfeucht-Fra\"\i ss\'e\ }

\documentclass[12pt]{article}
\usepackage{amssymb,latexsym}
\title{The Size of a Formula as a Measure of Complexity}
\author{
Lauri Hella \\ {\normalsize School of Information Sciences, Mathematics} \\ {\normalsize University of Tampere}
\and Jouko V\"a\"an\"anen\thanks{Research partially supported by
grant 251557 of the Academy of Finland and the EUROCORES LogICCC LINT programme.}\\ 
{\normalsize Department of Mathematics and Statistics} \\ {\normalsize University of Helsinki} \\ 
{\normalsize and} \\ {\normalsize Institute for Logic, Language and Computation} 
\\ {\normalsize University of Amsterdam}
}
\date{}
\begin{document}
\maketitle

\def\A{\mathcal{A}}
\def\B{\mathcal{B}}
\def\C{\mathcal{C}}
\def\D{\mathcal{D}}
\def\N{\mathcal{N}}
\def\O{\mathcal{O}}
\def\P{\mathcal{P}}
\def\S{{S}}
\def\R{{R}}
\def\K{\mathcal{K}}
\def\ma{\mathfrak{A}}
\def\map{\mathfrak{A'}}
\def\mb{\mathfrak{B}}
\def\mbp{\mathfrak{B'}}
\def\sma{{s}}
\def\smap{{s'}}
\def\smb{{r}}
\def\smbp{{r'}}
\def\ab{(\A,\B)}
\def\sab{(\S,\R)}
\newcommand{\abn}[1]{(\S_{#1},\R_{#1})}
\newcommand{\abpn}[1]{(\S'_{#1},\R'_{#1})}
\def\abp{(\A',\B')}
\newcommand{\sabn}[1]{(\S_{#1},\R_{#1})}
\newcommand{\sabpn}[1]{(\S'_{#1},\R'_{#1})}
\def\sabp{(\S',\R')}
\def\ba{(\B,\A)}
\newcommand{\ban}[1]{(\B_{#1},\A_{#1})}
\newcommand{\bapn}[1]{(\B'_{#1},\A'_{#1})}
\def\bap{(\B',\A')}
\def\sba{(\R,\S)}
\newcommand{\sban}[1]{(\R_{#1},\S_{#1})}
\newcommand{\sbapn}[1]{(\R'_{#1},\S'_{#1})}
\def\sbap{(\R',\S')}

\section{Introduction}

We propose a refinement of the usual Ehrenfeucht-Fra\"{\i}ss\' e game.
The new game will help us make finer distinctions than the
traditional one. In particular, it can be used to measure not only
quantifier rank but also lengths of conjunctions and disjunctions
needed for expressing a given property. Our game is similar to the
game in \cite{MR1995002} and in \cite{MR1914464}.

The most common measure of complexity of a first order sentence is
its quantifier rank, and the method of Ehrenfeucht-Fra\"{\i}ss\' e
games can be successfully used to prove 
lower 
bound results for
this measure of complexity. However, the number of non-equivalent
first order formulas of quantifier rank $n$ is an exponential
tower of height $n$
with a polynomial depending on the vocabulary
on the top. Thus we are not very close to knowing the formula if
we merely know its quantifier rank. In this paper we measure
complexity in terms of the number of symbols in the formula. The
advantage of our measure is that there are only $2^{P(n)}$
non-equivalent formulas with $n$ symbols, where $P(x)$ is a
polynomial depending on the vocabulary.

Let us consider  the question of the complexity of deciding
whether a binary string
\begin{equation}\label{string}s_1s_2\ldots s_n=10110010...01\end{equation} 
has a certain
pattern. One approach to this is circuit complexity. For example,
it is proved in \cite{MR738749} that constant depth circuits
cannot decide the {\em parity} of (\ref{string}), i.e. the
question whether (\ref{string}) has an even number of ones.
Another approach is computational complexity, where the parity of
(\ref{string}) is easily decided in linear time. In so called
descriptive complexity the question is posed, what is the logical
complexity of the simplest formula in a given logic that expresses
the property of (\ref{string}) in question, for example parity.
Concerning this type of problems, it was proved in \cite{MR706289}
that the question whether a set $R$ of
$n+1$-sequences\begin{equation}\label{relation}(a_1,\ldots,a_{n+1})
\end{equation} from a finite set $\{1,\ldots,m\}$
has an even cardinality cannot be expressed in  existential second
order logic, where the second order variables range over $\le
n$-ary relations, that is, there is no sentence $\phi$ of such
existential second order logic with the property that the
structure $(\{1,\ldots,m\},R)$ satisfies $\phi$ if and only if $R$
has even cardinality.

Let us return to (\ref{string}). A natural setup of using logic to
study the question of complexity of finding patterns in
(\ref{string}) is to use the propositional logic with
propositional symbols
\begin{equation}\label{psymb}p_1,p_2,\ldots,p_n\end{equation} with the meaning
\begin{equation}\label{prop}p_i\mbox{ is true if and only if
$s_i=1$.}
\end{equation}Propositional logic gives rise to
a variety of complexity measures. We focus here on the simplest of
them, namely the number of occurrences of propositional symbols in
a formula. Up to a constant factor this is the same as the number
of symbols, that is, the length of the formula.

 The parity of
(\ref{string}) is naturally expressed by the propositional formula
%
%
\begin{equation}\label{parity}\bigvee\{
\bigwedge_{i=1}^n p_i^A:A\subseteq\{1,\ldots,n\}, |A| \mbox{even}\},
\end{equation}
where 
\begin{equation}\label{ei}
p^A_i =
\cases{p_i & \mbox{if  $i\in A$}
\cr
\neg p_i  & \mbox{if  $i\notin A$,}}
\end{equation} 
which
gives the upper bound $\frac{1}{2}n2^n$ to the size of the
smallest formula expressing parity. However, there is a shorter
formula $\psi$ defined as follows: 
%
%
$\psi=\phi_{1,n}$, where for $1\le i<j\le n$ and $k=\lfloor\frac{i+j}{2}\rfloor$
\begin{eqnarray*}
\phi_{i,j}&=&(\phi_{i,k}\wedge
\phi_{k+1,j})\vee(\neg\phi_{i,k}\wedge
\neg\phi_{k+1,j})
\end{eqnarray*} and 
$$\phi_{i,i}=\neg p_i.$$ 
This gives (up to a constant
factor) the upper bound $n^2$ for the shortest propositional
formula expressing parity\footnote{More precisely, a straightforward 
induction shows that the size of $\s$ is at most $(n+1)^2$,
and if $n$ is a power of $2$, then the size of $\s$ is exactly $n^2$.}. 
Krapchenko (see \cite[p 258]{MR905473})
has proved  in 1971 that this is optimal. We prove below the same
result using our refinement of the Ehrenfeucht-Fra\"{\i}ss\' e
game.

For a slightly different complexity question, suppose we have n
unary predicates:
\begin{equation}\label{unary}P_1(x),\ldots,P_n(x).
\end{equation}

We ask what is the length of the shortest formula in first order
logic with the unary predicates (\ref{unary}) that expresses the
property  that every Boolean combination of the predicates is
non-empty. This property can be written as the existential
sentence:
%
%

\begin{equation}\label{exist}\bigwedge\{\exists x
\bigwedge_{i=1}^nP_i^A(x):A\subseteq\{1,\ldots,n\}\},
\end{equation}
where, as above in (\ref{ei}), $P_i^A(x)=P_i(x)$ if $i\in A$, and
$P^A_i(x)=\neg P_i(x)$ if $i\not\in A$. 

This gives the upper bound $(n+1)2^n$ to the
size of the smallest sentence expressing the given property
(we define the size of a formula to be the total number of atomic 
subformulas and quantifiers occurring in it). In Section 
\ref{exist}, we
show that this is optimal for existential sentences. If we give up
existential sentences and allow a quantifier alternation, we have a
smaller sentence $\phi$ defined as follows:
\begin{eqnarray*}
&&\forall x\exists y\bigwedge_{i=1}^n(P_i(x)\leftrightarrow
P_{i+1}(y))\wedge\\
&&\forall x\exists y(\bigwedge_{i=2}^n(P_i(x)\leftrightarrow
P_{i}(y))\wedge (P_1(x)\leftrightarrow \neg P_{1}(y))),
\end{eqnarray*}
where addition $i+1$ is to be calculated modulo $n$. The size of
$\phi$ is only $8n+4$, so this is really optimal 
(up to a constant factor) because we cannot
have such a $\phi$ without mentioning each $P_i$ at least once.


In our final application, we consider the minimal size of a 
first-order sentence expressing that the length of a linear order 
is at least $n$. It is well known that this can be expressed
by a sentence with quantifier rank $\lceil\log n\rceil+1$.
Such a sentence $\phi_n$ can be obtained as follows: define
first recursively formulas $\theta_k(x,y)$ 
by letting $\theta_2(x,y):= x<y$ and  
$\theta_k(x,y):= \exists z(\theta_l(x,z)\land\theta_m(z,y))$, where
$l=\lfloor k/2\rfloor +1$ and $m=\lceil k/2\rceil$, for $k>2$. Thus,
$\theta_k(x,y)$ says that there are at least $k$ elements in the interval $[x,y]$.
Then $\phi_n$ can be defined as $\exists x\exists y \theta_n(x,y)$.

Clearly the size of $\phi_n$ is linear in $n$. However, with a clever trick of recycling 
quantified variables it is possible to define a variant $\phi'_n$ of $\phi_n$ that has size 
$c \lceil\log n\rceil$, where $c$ is a small constant. Moreover, $\phi'_n$ can be 
chosen to be in the $4$-variable logic ${\rm FO}_4$ (see \cite{MR2274625}, Section 5).
On the other hand, Grohe and Schweikardt prove in \cite{MR2274625} that the minimal
size of an equivalent sentence in ${\rm FO}_3$ is at least ${1\over 2}\sqrt n$. Their proof
is based on a method which is related to the game we study in this paper, 
but it is not formulated in a game theoretic form.

\section{A Game for Propositional Logic}

We introduce now a game for measuring how long a propositional
formula has to be to express a given property of binary strings
\begin{equation}\label{strings}s=s_1\ldots s_n,\end{equation} where
$s_i\in\{0,1\}$. We denote strings by $s,r,t$ etc. Let $n$ be the
fixed length of the strings considered. We shall discuss nonempty
sets $\S$ of such strings and call them {\em string properties}.
For example, $S$ could be the set of strings (\ref{strings}) where
the cardinality of the set $\{i\in\{1,\ldots,n\}:s_i=1\}$ is even.
Another example is  the set of strings (\ref{strings}) with more
ones than zeros. In an extreme case  $S$ can be a singleton, e.g.
the singleton constant one string $\{111 \ldots 1\}$.

Propositional formulas are built up from propositional symbols
(\ref{psymb}) by means of $\neg, \wedge$ and $\vee$. We use Greek
letters $\phi,\psi$ etc to denote propositional formulas. The
concept $s\models\phi$ of a string $s$ satisfying a propositional
formula $\phi$ is defined in the usual way: $s\models p_i$ if
$s_i=1$; $s\models \neg\phi$ if $s\not\models \phi$; $s\models
\phi\wedge\psi$ if $s\models \phi$ and $s\models \psi$; $s\models
\phi\vee\psi$ if $s\models \phi$ or $s\models \psi$.

\begin{definition} Suppose $\phi$ is a
 propositional formula and $S$ a string property. We say that
  $\phi$ {\em defines}
 $S$ if for all strings $s$ $$s\in S\iff
s\models\phi.$$
\end{definition}

We define the {\em size} $\rk(\p)$ of the propositional formula
$\phi$ as follows:
\begin{eqnarray*}
\rk(p_i)&=&1\\
\rk(\neg\p)&=&\rk(\p)\\
\rk(\p{ {\vee} }\s)&=&\rk(\p)+\rk(\s)\\
\rk(\p{ {\wedge} }\s)&=&\rk(\p)+\rk(\s).
\end{eqnarray*}

\noindent Note that $\rk(\phi)$ is just the number of occurrences
of propositional symbols in $\phi$. Up to a constant factor, this
is the same as the number of symbols in $\phi$.

Note that while there are, up to logical equivalence, $2^{2^n}$
propositional formulas 
over the propositional symbols $p_1,\ldots,p_n$
altogether, there are for each $m$, up to
logical equivalence, only at most
$2^m\cdot(n+2)^{2m}$ propositional formulas $\phi$ with $\rk(\phi)\le m$.

A natural question now is:
\begin{quote}{\bf Question: }
Given a string property $S$, what is the size of the smallest
propositional formula that defines $S$.
\end{quote}

\noindent We shall  define a game for the study of this question.
In defining the game we follow mostly \cite{MR1995002} and to a
lesser degree \cite{MR1914464}. This game resembles the usual
\EF game for the \fl, but also the \EF game for the \ifl,
presented in \cite{MR1914464}.

In the usual \EF game we have two structures $\ma$ and $\mb$ and
the game is about their similarity. Player II maintains that the
structures are very similar. Player I maintains that there is a
difference. During the game player I picks elements from either
one of the models and player II tries to find similar elements in
the other model.

In the new game we do not have just two models but instead two
classes of models. Player II maintains that whatever features the
models in one class  have all in common, some member of the other
class also has. The game to be defined is able to detect very
minute differences in models and is in that sense finer than the
usual Ehrenfeucht game. In particular, this game can detect
differences in the size of conjunctions and disjunctions needed
for separating the classes of models.



\begin{definition}\label{ef-string}
Let $S$ and $R$ be string properties and $w$ a positive integer. The
game \(\ef_w\sab\) has two players. The number \(w\) is called the
{\em rank} of the game. In the beginning the position is
\((w,S,R)\). Suppose the position after $m$ moves is
\((w_m,S_m,R_m)\). There are the following two possibilities for
the continuation of the game:
\begin{description}

\item[Left splitting move:] Player I first chooses numbers
$u$ and $v$ such that $1\le u,v <w$ and $u+v=w$. Then
player I represents $\S_{m}$ as a union $C\cup D$. Now  the game
continues from the position \((u,C,R_m)\) or from the position
\((v,D,R_m)\), and player II can choose which.

\item[Right splitting move:] Player I first chooses numbers
$u$ and $v$ such that $1\le u,v <w$ and $u+v=w$. Then
Player I represents $\R_{m}$ as a union $C\cup D$. Now the game
continues from the position \((u,\S_m,C)\) or from the position
\((v,\S_m,D)\), and player II can choose which. 

\end{description}

\noindent The
game ends in a position \((w_m,S_m,R_m)\)
and player I wins if there is a propositional symbol $p_i$ such that 
either $s\models p_i$ for all $s\in S_m$ and $r\not\models p_i$ for all
$r\in R_m$, or $s\not\models p_i$ for all $s\in S_m$ and $r\models p_i$ 
for all $r\in R_m$. Player II wins the game if they reach a position 
\((w_m,S_m,R_m)\) such that $w_m=1$ and player I does not win in 
this position.
\end{definition}

This is  a game of perfect information and the concept of winning
strategy is defined as usual. Since the rank $w_m$ decreases in each
move, the game always ends in a finite number of moves.
Hence the game is determined: exactly one of the players has a winning 
strategy in $\ef_w(S,R)$. 

Note that if the players have reached a position $(w_m,S_m,R_m)$ 
in the game $\ef_w(S,R)$, then the continuation of the game
from that position onwards can be seen as a play of the
game $\ef_{w_m}(S_m,R_m)$. Thus, player I (player II)
has a winning strategy in the game $\ef_w(S,R)$ {\em 
from position} $(w_m,S_m,R_m)$ onwards if an only if player I (player II,
respectively) has a winning startegy in the game $\ef_{w_m}(S_m,R_m)$.

We say that a propositional formula $\p$ {\em separates} string properties
$S$ and $R$, in symbols $(S,R)\models\p$, if $s\models\p$ for all 
$s\in S$ and $r\not\models\p$ for all $r\in R$. 

\begin{theorem}\label{2}
  Suppose \(\sab\) is a pair of string properties, and let
$w$ be a positive integer.
Then the following conditions are equivalent:
\begin{description}

\item[\((1)_w\)] Player I has a winning strategy in the game
\(\ef_w\sab\).

\item[\((2)_w\)] There is a propositional formula $\p$ of size $\le
w$ such that $(S,R)\models\p$.
\end{description}

\end{theorem}

\proof 
We prove the equivalence of \((1)_w\) and \((2)_w\) by induction on $w$.
Consider first the case $w=1$.
By Definition~\ref{ef-string}, there are no moves in the game $\ef_1(S,R)$, 
and player I wins if and only if there is a proposition symbol
$p_i$ such that either $(S,R)\models p_i$, or $(S,R)\models \neg p_i$.
Since $\rk(p_i)=\rk(\neg p_i)=1$, we have \((1)_1\Longrightarrow (2)_1\).
On the other hand, if $\rk(\p)\le 1$, then $\p$ is, up to logical 
equivalence, either a proposition symbol, or the negation of a proposition
symbol. Thus, \((2)_1\) implies \((1)_1\).

Let us then consider the case \(w>1\), and assume \((1)_{v}\iff (2)_{v}\) 
for all \(v<w\) as
an induction hypothesis. To prove \((1)_{w}\iff (2)_{w}\), assume
first that player I has a winning strategy in the game
\(\ef_{w}(S,R)\). We have the following two cases according to the type of
the first move in the winning strategy of player I. 
\medskip

\noindent {\bf Case 1:} Player I makes a left splitting move by choosing 
$u$, $v$, $C$
and $D$ such that $1\le u,v<w$, $w=u+v$ and $S=C\cup D$. Since this move is
given by a winning strategy, player I has a winning strategy in both
of the games $\ef_u(C,R)$ and $\ef_v(D,R)$.
By induction hypothesis,
there are formulas $\s$ and $\h$ such that $\rk(\s)\le u$,
$\rk(\h)\le v$, $(C,R)\models\s$ and $(D,R)\models\h$. 
Thus, $s\models\s$ for all $s\in C$ and $s\models\h$ for all $s\in D$,
whence $s\models\s\lor\h$ for all $s\in S$. On the other hand, for all
$r\in R$, $r\not\models\s$ and $r\not\models\h$, whence consequently
$r\not\models\s\lor\h$. It follows that
$(S,R)\models\s\vee\h$. As $\rk(\s\vee\h)=\rk(\s)+\rk(\h)\le u+v=w$, 
$(2)_w$ is true.
\medskip

\noindent {\bf Case 2:} Player I makes a right splitting move by choosing  
$u$, $v$, $C$
and $D$ such that $1\le u,v<w$, $w=u+v$ and $R=C\cup D$. Since this move is
given by a winning strategy, player I has a winning strategy in both
of the games $\ef_u(S,C)$ and $\ef_v(S,D)$. By induction hypothesis
there are formulas $\s$ and $\h$ such that $\rk(\s)\le u$,
$\rk(\h)\le v$, $(S,C)\models\s$ and $(S,D)\models\h$. 
Thus, for all
$s\in S$, $s\models\s$ and $s\models\h$, whence 
$s\models\s\land\h$. On the other hand, 
$r\not\models\s$ for all $r\in C$ and $r\not\models\h$ for all $r\in D$,
whence $r\not\models\s\land\h$ for all $r\in R$.
It follows that
$(S,R)\models\s\wedge\h$. As $\rk(\s\wedge\h)=\rk(\s)+\rk(\h)\le u+v=w$, 
$(2)_w$
is true.
\medskip

To prove the converse implication
\((2)_{w}\Longrightarrow (1)_{w}\), assume that there is a formula
$\p$ of size $\le w$ such that $(S,R)\models\p$. We show that then
player I has a winning strategy in the game \(\ef_{w}(S,R)\). 
We assume without loss of generality that $\p$ is in {\em negation
normal form}\footnote{A formula is in negation normal form, if
all negations occurring in it are in front of proposition symbols.
A straightforward induction shows that for any formula $\p$
there is an equivalent formula $\p'$ in negation normal form 
such that $\rk(\p')=\rk(\p)$.}. The
strategy depends on $\p$ as follows.
\medskip

\noindent {\bf Case 1:} $\p$ is a proposition symbol or the negation of
a proposition symbol. Then by the definition of the game $\ef_w$,
player I wins without making any moves.
\medskip

\noindent {\bf Case 2:} $\p$ is $\s\vee\h$.
Let $C$ be the set of strings
$s\in S$ such that $s\models\s$, and let $D$ be the set of strings
$s\in S$ such that $s\models\h$. Since $s\models\p$ for every  
$s\in S$, we have $S=C\cup D$. Moreover, since $r\not\models\p$ 
for every 
$r\in R$, we have $(C,R)\models\s$ and $(D,R)\models\h$. Finally, as 
$\rk(\p)\le w$, there are $u$ and $v$ such that $w=u+v$, $\rk(\s)\le u$
and $\rk(\h)\le v$. Note that $u,v<w$, and so, by induction hypothesis,
player I has a winning strategy in the games \(\ef_{u}(C,R)\) and
\(\ef_{v}(D,R)\). Thus, using $w=u+v$, $S=C\cup D$ as his first move
and these winning strategies in the rest of the game,
player I is guaranteed to win the game \(\ef_{w}(S,R)\).
\medskip

\noindent {\bf Case 3:} $\p$ is $\s\wedge\h$.
Let $C$ be the set of strings
$r\in R$ such that $r\not\models\s$, and let $D$ be the set of strings
$r\in R$ such that $s\not\models\h$. Since $r\not\models\p$ for every  
$r\in R$, we have $R=C\cup D$. Moreover, since $s\models\p$ for every  
$s\in S$, we have $(S,C)\models\s$ and $(S,D)\models\h$. Finally, as 
$\rk(\p)\le w$, there are $u$ and $v$ such that $w=u+v$, $\rk(\s)\le u$
and $\rk(\h)\le v$. Note that $u,v<w$, and so, by induction hypothesis,
player I has a winning strategy in the games \(\ef_{u}(S,C)\) and
\(\ef_{v}(S,D)\). Thus, using $w=u+v$, $R=C\cup D$ in his first move,
and these winning strategies in the rest of the game,
player I is guaranteed to win the game \(\ef_{w}(S,R)\).

\qed

\begin{corollary}
 Suppose $K$ is a string property.
Then the following conditions are equivalent:
\begin{description}
\item[(1)]  \(K\) is definable by
 a propositional formula $\p$ of size $\le w$.

\item[(2)] Player I has a winning strategy in the game
\(\ef_w\sab\) for all $\S$ and $\R$ such that $\S\subseteq K$ and
$\R\cap  K=\emptyset$.
\end{description}
 \end{corollary}

\proof Suppose $K$ is defined by $\p$ of size $\le w$. Let $\S$
and $\R$ be given so that $\S\subseteq K$ and
$\R\cap\K=\emptyset$. Thus every model in $\S$ satisfies $\p$ and
no model in $\R$ does, whence $\sab\models\p$. By
Theorem~\ref{2}, player I has a winning strategy in $\ef_w\sab$.
For the converse, suppose player I has a winning strategy in
$\ef_w\sab$ for all $\sab$ with $\S\subseteq\K$ and $\R\cap K$.
This holds in particular if $\S$ is $K$ and $\R$ is the complement of $K$. 
Thus for this
choice of $\S$ and $\R$ we have $\sab\models\p$ for some $\p$ of
size $\le w$. Thus $K$ is defined by $\p$. \qed

We get the following method for showing that a string property $K$
is {\em not} definable by a propositional formula of size $\le w$:
We find  classes $\S$ and $\R$ such that
\begin{enumerate}

\item $\S\subseteq K$.

\item $\R\cap K=\emptyset$.

\item Player II has a winning strategy in $\ef_w\sab$.

\end{enumerate}

\section{The Propositional Complexity of Parity}

We shall now prove that the parity of a binary string of $n$ bits
cannot be expressed with a propositional formula smaller than
$n^2$. The result has been proved in \cite{MR905473} with a method
which is very similar to ours, even if it does not explicitly use
games. We present the proof as an example of the use of our game.

Suppose $S$ and $R$ are disjoint sets of binary strings of length
$n$. Let $$E=\{(f,g)\in S\times R:|\{i:f_i\ne g_i\}|=1\}.$$ The
{\em density} of the pair $(S,R)$ is defined to be the pair
$$D(S,R)=(s,r),$$ where $$s=\frac{|E|}{|S|}, \mbox{ and }
r=\frac{|E|}{|R|}.$$

\begin{lemma}\label{base}
If $s>1$ or $r>1$, then there is no proposition symbol $p_i$
such that $(S,R)\models\ p_i$ or $(S,R)\models\neg p_i$.
\end{lemma}

\proof Suppose $s>1$ and $i\in\{1,\ldots,n\}$. Then there are $f\in S$ 
and $g,h\in R$ such
that $(f,g),(f,h)\in E$. Clearly this means that $f_i=g_i$ or $f_i=h_i$,
whence $f\models p_i\iff g\models p_i$ or 
$f\models p_i\iff h\models p_i$. Thus, if all strings in $S$
satisfy $p_i$ ($\neg p_i$), then there is a string in $R$ 
satisfying $p_i$ ($\neg p_i$, respectively). It follows that neither $p_i$,
nor $\neg p_i$ separates the sets $S$ and $R$.
 
In the case $r>1$ there are $f,g\in S$ and $h\in R$ such that
$(f,h),(g,h)\in E$. Then we have $h_i=f_i$ or $h_i=g_i$,
whence $h\models p_i\iff f\models p_i$ or 
$h\models p_i\iff g\models p_i$. As above, we conclude that neither $p_i$,
nor $\neg p_i$ separates the sets $S$ and $R$.
\qed

\begin{lemma}\label{str}
(a) Suppose $D(S,R)=(s,r)$ and $S=S_0\cup S_1$, where $S_0\cap S_1=
\emptyset$. Let
$D(S_0,R)=(s_0,r_0)$ and $D(S_1,R)=(s_1,r_1)$. Then
$s_0r_0+s_1r_1\ge sr$.

(b) Suppose $D(S,R)=(s,r)$ and $R=R_0\cup R_1$, where $R_0\cap R_1=
\emptyset$. Let
$D(S,R_0)=(s_0,r_0)$ and $D(S,R_1)=(s_1,r_1)$. Then
$s_0r_0+s_1r_1\ge sr$.
\end{lemma}

\proof We prove (a); the proof of (b) is similar. Clearly $r_0+r_1=r$ and
$\frac{r_0}{s_0}+\frac{r_1}{s_1}=\frac{r}{s}$. Since the harmonic
mean is never greater than the arithmetic mean\footnote{We are
indebted to Kerkko Luosto for pointing this out.}, we have 
$$\frac{s_0r_0+s_1r_1}{r}=\frac{s_0r_0+s_1r_1}{r_0+r_1}\ge
\frac{r_0+r_1}{\frac{r_0}{s_0}+\frac{r_1}{s_1}}=
\frac{r}{\frac{r}{s}}
=s,$$
and the claim
follows. \qed

\begin{lemma}\label{ws}
Suppose $D(S,R)=(s,r)$ and $w<sr$. Then player II has a 
winning strategy in $\ef_w(S,R)$. 
\end{lemma}

\proof The proof is by induction on $w$. In the case $w=1$,
the assumption $w<sr$ implies that either $s>1$ or
$r>1$. By Lemma~\ref{base}, there is no proposition symbol $p_i$
such that $(S,R)\models p_i$ or $(S,R)\models\neg p_i$. 
Since the game $\ef_1(S,R)$ ends in its initial position 
$(1,S,R)$, this means that player II automatically wins 
this game.

Assume then that $w>1$. Using Lemma~\ref{base} again, we see that player I
does not win the game $\ef_w(S,R)$ without making moves.
Suppose then that player I makes a left splitting move 
$w=u+v$ and $S=C\cup D$.
Let $S_0\subseteq C$ and $S_1\subseteq D$ be sets such that $S=S_0\cup S_1$
and $S_0\cap S_1=\emptyset$.
Let $D(S_0,R)=(s_0,r_0)$ and $D(S_1,R)=(s_1,r_1)$. Then by
Lemma~\ref{str}(a), $s_0r_0+s_1r_1\ge sr$. Thus $s_0r_0> u$ or
$s_1r_1> v$, for otherwise $$sr\le s_0r_0+s_1r_1\le u+v=w,$$
contrary to the assumption. If $s_0r_0>u$, then
by induction hypothesis, player II
has a winning strategy in the game $\ef_{u}(S_0,R)$, whence,
a fortiori, he has a winning strategy in $\ef_{u}(C,R)$.
Otherwise, $s_1r_1> v$, and by the same argument, player II
has a winning strategy in the game $\ef_{v}(D,R)$. Thus, in 
any case, player II can make his move in such a way that he is guaranteed
to win. 

The case of a right splitting move $w=u+v$ and $R=C\cup D$
is proved in the same way by using Lemma~\ref{str}(b). \qed

We are now ready to prove the promised lower bound for parity of
binary strings.

\begin{corollary} [\cite{MR905473}]
If $\p$ is a propositional formula expressing the parity 
of strings $s\in\{0,1\}^n$, then the size of $\p$ is at least $n^2$. 
\end{corollary}

\proof Let $S$ be the set of all strings $f\in\{0,1\}^n$ such that 
$|\{i:f_i=1\}|$ is even, and let $R$ be the complement of $S$.
Thus, $|S|=|R|=2^{n-1}$. Furthemore, for each string $f\in S$,
there are $n$ different strings $g\in\{0,1\}^n$ such that 
$|\{i:f_i\not=g_i\}|=1$, and all these strings are in $R$. Thus, we have
$|E|=2^{n-1}n$, and consequently $D(S,R)=(n,n)$. By Lemma~\ref{ws},
player II has a winning strategy in the game $\ef_w(S,R)$ whenever
$w<n^2$, and we conclude that $S$ is not definable by any formula
of size less than $n^2$.
\qed

As noted in the introduction, the parity of strings in $\{0,1\}^n$
can be expressed by a formula of size at most $(n+1)^2$. Thus, the 
lower bound $n^2$ cannot be essentially improved.

\section{A Game for Predicate Logic}

We shall next define a game that can be used for measuring
the size of a first-order sentence needed for expressing properties
of models. For the sake of simplicity, we will only consider models with
relational vocabulary.

We need to fix some notation first. The
universe of a model $\ma$ is denoted by $A$, of $\mb$ by $B$, etc.
We use $x_j$, $j\in\oN$, to denote variables. 
A variable assignment for a model $\ma$ is a finite partial mapping
$\alpha:\oN\to A$. The finite domain of $\alpha$ is denoted by
$\dom(\alpha)$. If $\p$ is a formula, then $(\ma,\alpha)\models\p$ means
that the assignment $\alpha$ satisfies the formula $\p$ in the
model $\ma$. For this to be meaningful the domain of the
assignment $\alpha$ has to include all the $j$ for which the variable 
$x_j$ is free in $\p$.

We shall discuss classes $\A$ of structures $(\ma,\alpha)$,
where $\ma$ is a model and $\alpha$ is an assignment. 
We assume that whenever $(\ma,\alpha),(\mb,\beta)\in\A$, then $\ma$ and
$\mb$ have the same vocabulary, and $\alpha$ and $\beta$ have
the same domain, which we denote by $\dom(\A)$. If $\alpha$ is an
assignment on $\ma$,
$a\in A$ and $j\in\oN$, then $\alpha(a/j)$ is the
assignment that maps $j$ to $a$ and agrees
with $\alpha$ otherwise.
If $F$ is a choice function on $\A$,
that is, $F$ is a function defined on $\A$ such that 
$F((\ma,\alpha))\in A$ for all $(\ma,\alpha)\in\A$, then $\A(F/j)$ is
defined as $\{(\ma,\alpha(F((\ma,\alpha))/j):(\ma,\alpha)\in\A\}$. Finally,
$\A(\star/j)=\{(\ma,\alpha(a/j)):(\ma,\alpha)\in\A, a\in A\}$.

Let $\A$ and $\B$ be classes of structures of a fixed relational vocabulary.
Assume further that $\dom(\A)=\dom(\B)$, and $\p$ is a formula such that
$j\in\dom(\A)$ for all variables $x_j$ which are free in $\p$.
As in the case of propositional logic, we say that $\p$ {\em separates}
the classes $\A$ and $\B$, $(\A,\B)\models\p$, if $(\ma,\alpha)\models\p$
for all $(\ma,\alpha)\in\A$ and $(\mb,\beta)\not\models\p$
for all $(\mb,\beta)\in\B$.

\begin{definition}\label{foef}
Let $\A$ and $\B$ be classes of structures of the same relational 
vocabulary with $\dom(\A)=\dom(\B)$, and let $w$ be a 
positive integer.
The game \(\ef_w\ab\) has two
players. The number \(w\) is called the {\em rank} of the game. In
the beginning the position is \((w,\A,\B)\). Suppose the position
after $m$ moves is \((w_m,\A_m,\B_m)\), where
$\dom(\A_m)=\dom(\B_m)$. There are the following four
possibilities for the continuation of the game:
\begin{description}

\item[Left splitting move:] Player I first chooses numbers
$u$ and $v$ such that $1\le u,v<w$ and $u+v=w_m$. Then
Player I represents $\A_{m}$ as a union $\C\cup \D$. Now the game
continues from the position \((u,\C,\B_m)\) or from the position
\((v,\D,\B_m)\), and player II can choose which. 

\item[Right splitting move:] Player I first chooses numbers
$u$ and $v$ such that $1\le u,v<w$ and $u+v=w_m$. Then
Player I represents $\B_{m}$ as a union $\C\cup \D$. Now the game
continues from the position \((u,\A_m,\C)\) or from the position
\((v,\A_m,\D)\), and player II can choose which. 

\item[Left supplementing move:] Player I chooses a natural number $j$
and a choice
function \(F\) for $\A_m$.  Then the game continues from the
position \((w_m-1,\A_m(F/j),\B_m(\star/j))\).

\item[Right supplementing move:] Player I chooses a natural number $j$
and a choice
function \(F\) for $\B_m$.  Then the game continues from the
position \((w_m-1,\A_m(\star/j),\B_m(F/j))\).

\end{description}

\noindent The game ends in a position \((w_m,\A_m,\B_m)\)
and player I wins if there is an atomic or a negated atomic
formula $\p$ such that 
$(\A_m,\B_m)\models\p$. Player II wins the game if they reach a position 
\((w_m,S_m,R_m)\) such that $w_m=1$ and player I does not win in 
this position.
\end{definition}

This is  a game of perfect information and the concept of winning
strategy is defined as usual. The game is determined by the
Gale-Stewart theorem.

We define the {\em size} $\rk(\p)$ of the formula $\phi$ of
predicate logic as follows:

\begin{eqnarray*}
\rk(\p)&=&1\;\mbox{ for atomic $\p$}\\
\rk(\neg\p)&=&\rk(\p)\\
\rk(\p{ {\vee} }\s)&=&\rk(\p)+\rk(\s)\\
\rk(\p{ {\wedge} }\s)&=&\rk(\p)+\rk(\s)\\
\rk(\exists x_j\p)&=&\rk(\p)+1\\
\rk(\forall x_j\p)&=&\rk(\p)+1\end{eqnarray*}

Note that there are for each $w$, up to logical
equivalence, only finitely many formulas of size $\le w$.

\begin{theorem}\label{2p}
  Suppose \(\ab\) is a pair of classes
  of structures of the same vocabulary,
and let $w$ be a positive integer.
Then the following conditions are equivalent:
\begin{description}

\item[\((1)_w\)] Player I has a winning strategy in the game
\(\ef_w\ab\).

\item[\((2)_w\)] There is a formula $\p$ of predicate logic of size
$\le w$ such that $(\A,\B)\models\p$.

\end{description}

\end{theorem}

\proof 
We prove the equivalence of \((1)_w\) and \((2)_w\) by induction
on $w$. Consider first the case \(w=1\). By Definition~\ref{foef},
there are no moves in the game $\ef_1(\A,\B)$, and player I wins
if and only if there is an atomic formula $\p$ such that either
$(\A,\B)\models\p$ or $(\A,\B)\models\neg\p$. Since $\rk(\p)=
\rk(\neg\p)=1$, we have \((1)_1\Longrightarrow (2)_1\). On the other
hand, if $\rk(\p)=1$, then $\p$ is, up to logical equivalence,
either an atomic formula, or the negation of an atomic formula.
Thus, \((2)_1\) implies \((1)_1\).

Let us then consider the case $w>1$, and assume \((1)_v\iff (2)_v\) 
for all $v<w$ as an
induction hypothesis. To prove \((1)_w\iff (2)_w\), assume first that
player I has a winning strategy in the game $\ef_w(\A,\B)$. We have
the following four cases according to the type of the first move
in the winning strategy of player I. 
\medskip

\noindent {\bf Case 1:} Player I makes a left splitting move by choosing 
$u$, $v$, $\C$
and $\D$ such that $1\le u,v<w$, $w=u+v$ and $\A=\C\cup \D$. 
Since this move is
given by a winning strategy, player I has a winning startegy in both
of the games $\ef_u(\C,\B)$ and $\ef_v(\D,\B)$.
By induction hypothesis,
there are formulas $\s$ and $\h$ such that $\rk(\s)\le u$,
$\rk(\h)\le v$, $(\C,\B)\models\s$ and $(\D,\B)\models\h$. 
Thus, $(\ma,\alpha)\models\s$ for all $(\ma,\alpha)\in \C$ 
and $(\ma,\alpha)\models\h$ for all $(\ma,\alpha)\in \D$,
whence $(\ma,\alpha)\models\s\lor\h$ for all $(\ma,\alpha)\in \A$. 
On the other hand, for all
$(\mb,\beta)\in \B$, $(\mb,\beta)\not\models\s$ and 
$(\mb,\beta)\not\models\h$, whence consequently
$(\mb,\beta)\not\models\s\lor\h$. It follows that
$(\A,\B)\models\s\vee\h$. As $\rk(\s\vee\h)=\rk(\s)+\rk(\h)\le u+v=w$, 
$(2)_w$ is true.
\medskip

\noindent {\bf Case 2:} Player I makes a right splitting move by choosing  
$u$, $v$, $\C$
and $\D$ such that $1\le u,v<w$, $w=u+v$ and $\B=\C\cup \D$. 
Since this move is
given by a winning strategy, player I has a winning startegy in both
of the games $\ef_u(\A,\C)$ and $\ef_v(\A,\D)$. By induction hypothesis
there are formulas $\s$ and $\h$ such that $\rk(\s)\le u$,
$\rk(\h)\le v$, $(\A,\C)\models\s$ and $(\A,\D)\models\h$. 
Thus, for all
$(\ma,\alpha)\in \A$, $(\ma,\alpha)\models\s$ and $(\ma,\alpha)\models\h$, 
whence $(\ma,\alpha)\models\s\land\h$. On the other hand, 
$(\mb,\beta)\not\models\s$ for all $(\mb,\beta)\in \C$ 
and $(\mb,\beta)\not\models\h$ for all $(\mb,\beta)\in \D$,
whence $(\mb,\beta)\not\models\s\land\h$ for all $(\mb,\beta)\in \B$.
It follows that
$(\A,\B)\models\s\wedge\h$. As $\rk(\s\wedge\h)=\rk(\s)+\rk(\h)\le u+v=w$, 
$(2)_w$
is true.
\medskip

\noindent {\bf Case 3:} Player I makes a left supplementing move by choosing  
a natural number $j$ and a choice function $F$ for $\A$. The next 
position in the game is then $(w-1,\A(F/j),\B(\star/j))$.
Since this move is
given by a winning strategy, player I has a winning strategy in 
the game $\ef_{w-1}(\A(F/j),\B(\star/j))$. By induction hypothesis
there is a formula $\s$ such that $\rk(\s)\le w-1$ and 
$(\A(F/j),\B(\star/j))\models\s$. Let $\p$ be the formula
$\exists x_j\s$. Then $\rk(\p)=\rk(\s)+1\le w$, and it suffices
to show that $(\A,\B)\models\p$.
Note first that for all
$(\ma,\alpha)\in \A$, $(\ma,\alpha(a/j))\models\s$, where 
$a=F((\ma,\alpha))$. Thus we have $(\ma,\alpha)\models\p$ for
all $(\ma,\alpha)\in \A$.
On the other hand, for all $(\mb,\beta)\in \B$ and all $b\in B$, we have
$(\mb,\beta(b/j))\not\models\s$,
whence $(\mb,\beta)\not\models\p$. 
\medskip

\noindent {\bf Case 4:} Player I makes a right supplementing move by choosing  
a natural number $j$ and a choice function $F$ for $\B$. The next 
position in the game is then $(w-1,\A(\star/j),\B(F/j))$.
Since this move is
given by a winning strategy, player I has a winning strategy in 
the game $\ef_{w-1}(\A(\star/j),\B(F/j))$. By induction hypothesis
there is a formula $\s$ such that $\rk(\s)\le w-1$ and 
$(\A(\star/j),\B(F/j))\models\s$. Let $\p$ be the formula
$\forall x_j\s$. Then $\rk(\p)=\rk(\s)+1\le w$, and it suffices
to show that $(\A,\B)\models\p$.
Note first that for all
$(\ma,\alpha)\in \A$ and all $a\in A$, we have
$(\ma,\alpha(a/j))\models\s$. Thus we have $(\ma,\alpha)\models\p$ for
all $(\ma,\alpha)\in \A$.
On the other hand, for all $(\mb,\beta)\in \B$, 
$(\mb,\beta(b/j))\not\models\s$, where $b=F((\mb,\beta))$. Thus we have
$(\mb,\beta)\not\models\p$ for all $(\mb,\beta)\in \B$. 
\medskip

To prove the converse implication
\((2)_{w}\Longrightarrow (1)_{w}\), assume that there is a formula
$\p$ of size $\le w$ such that $(\A,\B)\models\p$. We show that then
player I has a winning strategy in the game \(\ef_{w}(\A,\B)\). 
As in the case of propositional logic, we can assume without 
loss of generality that $\p$ is in  negation
normal form. The
strategy of player I depends on $\p$ as follows.
\medskip

\noindent {\bf Case 1:} $\p$ is an atomic formula or the negation of
an atomic formula. Then by Definition \ref{foef}, 
player I wins the game $\ef_w(\A,\B)$ without making any moves.

\medskip

\noindent {\bf Case 2:} $\p$ is $\s\vee\h$.
Let $\C$ be the class of structures
$(\ma,\alpha)\in \A$ such that $(\ma,\alpha)\models\s$, and let $\D$ 
be the class of structures
$(\ma,\alpha)\in \A$ such that $(\ma,\alpha)\models\h$. 
Since $(\ma,\alpha)\models\p$ for every  
$(\ma,\alpha)\in\A$, we have $\A=\C\cup \D$. Moreover, 
since $(\mb,\beta)\not\models\p$ 
for every 
$(\mb,\beta)\in \B$, we have $(\C,\B)\models\s$ and $(\D,\B)\models\h$. 
Finally, as 
$\rk(\p)\le w$, there are $u$ and $v$ such that $w=u+v$, $\rk(\s)\le u$
and $\rk(\h)\le v$. Note that $u,v<w$, and so, by induction hypothesis,
player I has a winning strategy in the games \(\ef_{u}(\C,\B)\) and
\(\ef_{v}(\D,\B)\). Thus, using $w=u+v$, $\A=\C\cup \D$ as his first move
and these winning strategies in the rest of the game,
player I is guaranteed to win the game \(\ef_{w}(\A,\B)\).
\medskip

\noindent {\bf Case 3:} $\p$ is $\s\wedge\h$.
Let $\C$ be the class of structures
$(\mb,\beta)\in \B$ such that $(\mb,\beta)\not\models\s$, and let 
$\D$ be the class of structures
$(\mb,\beta)\in \B$ such that $(\mb,\beta)\not\models\h$. 
Since $(\mb,\beta)\not\models\p$ for every  
$(\mb,\beta)\in \B$, we have $\B=\C\cup \D$. Moreover, 
since $(\ma,\alpha)\models\p$ for every  
$(\ma,\alpha)\in \A$, we have $(\A,\C)\models\s$ and $(\A,\D)\models\h$. 
Finally, as 
$\rk(\p)\le w$, there are $u$ and $v$ such that $w=u+v$, $\rk(\s)\le u$
and $\rk(\h)\le v$. Note that $u,v<w$, and so, by induction hypothesis,
player I has a winning strategy in the games \(\ef_{u}(\A,\C)\) and
\(\ef_{v}(\A,\D)\). Thus, using $w=u+v$, $\A=\C\cup \D$ in his first move,
and these winning strategies in the rest of the game,
player I is guaranteed to win the game \(\ef_{w}(\A,\B)\).
\medskip

\noindent {\bf Case 4:} $\p$ is $\exists x_j\s$.
Since $(\ma,\alpha)\models\p$ for every $(\ma,\alpha)\in\A$, there
is a choice function $F$ for $\A$ such that 
$(\ma,\alpha(F((\ma,\alpha))/j))\models\s$ for all $(\ma,\alpha)\in \A$.
Thus, $(\ma,\alpha^*)\models\s$ for every $(\ma,\alpha^*)\in \A(F/j)$.
On the other hand, for all $(\mb,\beta)\in\B$ we have $(\mb,\beta)
\not\models\p$, whence $(\mb,\beta(b/j))\not\models\s$ for all $b\in B$.
In other words, $(\mb,\beta^*)\not\models\s$ for all $(\mb,\beta^*)\in
\B(\star/j)$. Thus we conclude that $(\A(F/j),\B(\star/j))\models\s$.
Note that $\rk(\s)=\rk(\p)-1\le w-1$, and so, by induction hypothesis,
player I has a winning strategy in the game $\ef_{w-1}(\A(F/j),\B(\star/j))$.
Thus, starting with the left supplementing move $j$ and $F$, and using
this winning strategy in the rest of the game, player I is guaranteed
to win the game $\ef_w(\A,\B)$.
\medskip

\noindent {\bf Case 5:} $\p$ is $\forall x_j\s$.
Since $(\mb,\beta)\not\models\p$ for every $(\mb,\beta)\in\B$, there
is a choice function $F$ for $\B$ such that 
$(\mb,\beta(F((\mb,\beta))/j))\not\models\s$ for all $(\mb,\beta)\in \B$.
Thus, $(\mb,\beta^*)\not\models\s$ for every $(\mb,\beta^*)\in \B(F/j)$.
On the other hand, for all $(\ma,\alpha)\in\A$ we have $(\ma,\alpha)
\models\p$, whence $(\ma,\alpha(a/j))\models\s$ for all $a\in A$.
In other words, $(\ma,\alpha^*)\models\s$ for all $(\ma,\alpha^*)\in
\A(\star/j)$. Thus we conclude that $(\A(\star/j),\B(F/j))\models\s$.
Note that $\rk(\s)=\rk(\p)-1\le w-1$, and so, by induction hypothesis,
player I has a winning strategy in the game $\ef_{w-1}(\A(\star/j),\B(F/j))$.
Thus, starting with the right supplementing move $j$ and $F$, and using
this winning strategy in the rest of the game, player I is guaranteed
to win the game $\ef_w(\A,\B)$.
 \qed

\begin{corollary}
 Suppose $\K$ is a class of models of the same vocabulary.
Then the following conditions are equivalent:
\begin{description}
\item[(1)]  \(\K\) is the class of models of a first order
sentence $\p$ of size $\le w$.

\item[(2)] Player I has a winning strategy in the game
\(\ef_w\ab\) for all $\A$ and $\B$ such that $\A\subseteq\K$ and
$\B\cap \K=\emptyset$.
\end{description}
 \end{corollary}

\proof Suppose $\K$ is the class of models of $\p$ and $\rk(\p)\le
w$. Let $\A$ and $\B$ be given so that $\A\subseteq \K$ and
$\B\cap\K=\emptyset$. Thus every model in $\A$ satisfies $\p$ and
no model in $\B$ does, whence $\ab\models\p$. By
Theorem~\ref{2p}, player I has a winning strategy in $\ef_w\ab$.
For the converse, suppose player I has a winning strategy in
$\ef_w\ab$ for all $\ab$ with $\A\subseteq\K$ and $\B\cap\K$. This
holds in particular if $\A=\K$ and $\B$ is the complement of $\K$. 
Thus for this choice
of $\A$ and $\B$ we have $\ab\models\p$ for some $\p$ of size $\le
w$. Thus $\K$ is the class of all models of $\p$. \qed

We get the following method for showing that a model class $\K$ is
{\em not} definable by a first order sentence of size  $\le w$: We
find classes $\A$ and $\B$ such that
\begin{enumerate}

\item $\A\subseteq \K$.

\item $\B\cap \K=\emptyset$.

\item Player II has a winning strategy in $\ef_w\ab$.

\end{enumerate}

\subsection*{A Game for Existential Formulas}

The game $\ef_w$ can be used for solving questions of the type 
\begin{quote}
Given a property of $\cal P$ models, what 
is the size of the smallest first order formula that defines $\cal P$?
\end{quote}
But in some applications we are interested in definability by some
restricted type of of formulas, rather than arbitrary first order
formulas. For example, if a property $\cal P$ of models is known
to be definable by an existential formula, it is natural
to ask, what is the size of the smallest existential formula defining 
$\cal P$. We will now define a variant of the game $\ef_w$ that
can be used in studying this question. 
 
Here we say that a first order formula is 
{\em existential} if it is in negation normal form, and it does not
contain any universal quantifiers. In other words,
existential formulas are built from atomic formulas and
negations of atomic formulas by using the connectives $\lor$, $\land$ and
the quantifier $\exists$.

Intuitively, in the definition of the game $\ef_w$, 
left and right splitting moves correspond to the
connectives $\lor$ and $\land$, while left and right
supplementing moves correspond to the quantifiers $\exists$
and $\forall$. Thus, we obtain a game for existential formulas
simply by dropping right supplementing moves.

\begin{definition}\label{eef}
Let $\A$ and $\B$ be classes of structures of the same relational 
vocabulary with $\dom(\A)=\dom(\B)$, and let $w$ be a 
positive integer.
The existential game \(\ef^\exists_w\ab\) has the same rules 
as $\ef_w(\A,\B)$, except that player I is not allowed to make
right supplementing moves.
\end{definition}

\begin{theorem}\label{egc}
  Suppose \(\ab\) is a pair of classes
  of structures of the same vocabulary,
and let $w$ be a positive integer.
Then the following conditions are equivalent:
\begin{description}

\item[\((1)_w\)] Player I has a winning strategy in the game
\(\ef^\exists_w\ab\).

\item[\((2)_w\)] There is an existential formula $\p$ of size
$\le w$ such that $(\A,\B)\models\p$.

\end{description}

\end{theorem}

\proof The claim is proved by a straightforward adaptation
of the proof of Lemma~\ref{2p}. In the proof of 
\((1)_w\Longrightarrow (2)_w\) it suffices to note that
omitting Case~4, the formula $\p$ separating the classes
$\A$ and $\B$ will always be existential. Similarly, in the
proof of the opposite implication, Case~5 never occurs as the
separating formula $\p$ is existential, whence the 
winning strategy of player I
does not use right supplementing moves.
\qed

As in the case of the full first order game $\ef_w$, we get again 
a method for showing that a model class $\K$ is
{\em not} definable by an existential formula of size  $\le w$: We
find classes $\A$ and $\B$ such that
\begin{enumerate}

\item $\A\subseteq \K$.

\item $\B\cap \K=\emptyset$.

\item Player II has a winning strategy in $\ef^\exists_w\ab$.

\end{enumerate}

\section{The Existential Complexity
of Non-Empti\-ness of Boolean Combinations}\label{exist}

We observed in the introduction that non-emptiness of all Boolean
combinations of $n$ unary predicates can be defined with a
sentence of size $8n+4$ if a quantifier alternation is allowed. If
only existential sentences are allowed, this can be expressed with
a sentence of size $(n+1)2^n$. We use now our game for existential
formulas to show that $(n+1)2^n$ is the best possible value.

For each binary string $s\in \{0,1\}^n$, let $\mb_s$ be the 
$\{P_1,\ldots,P_n\}$-structure such that for each $r\in\{0,1\}^n$, the 
Boolean combination
of $P_1^{\mb_s},\ldots,P_n^{\mb_s}$ corresponding to $r$ contains exactly
two elements, $b_r$ and $c_r$, except that the combination corresponding 
to $s$ is empty. Furthermore, let $\ma$ be the structure in which each
Boolean combination contains exactly one element; let $a_r$ be the
element in the Boolean combination corresponding to $r\in\{0,1\}^n$.  
Let $\A_0$ be the
class $\{(\ma,\emptyset)\}$, and let $\B_0$ be the class 
$\{(\mb_s,\emptyset): s\in \{0,1\}^n\}$.
We will show that player II has a winning strategy in the existential game
$\ef^\exists_w(\A_0,\B_0)$ for all $w<(n+1)2^n$.

Since the game is existential, in any position $(u,\A,\B)$, the set
$\A$ contains only one structure $(\ma,\alpha)$, where  $\alpha$
is the variable assignment that results from the moves of player I. 
We say that a structure
$(\mb_s,\beta)$ is {\em flawless} (with respect to $\alpha$), if 
$\dom(\beta)=\dom(\alpha)$, and for all 
$j\in\dom(\alpha)$
and all $r\in\{0,1\}^n$, we have 
\begin{equation}\label{ehto1}
\alpha(j)=a_r\iff\beta(j)=b_r.
\end{equation} 
Note that if $a_s\in\ran(\alpha)$, then there is no $\beta$
such that $(\mb_s,\beta)$ is flawless. On the other hand,
if $a_s\not\in\ran(\alpha)$, then there is a unique $\beta$
such that $(\mb_s,\beta)$ is flawless; we denote this $\beta$
by $\beta_{s,\alpha}$.

Furthermore, we say that a structure
$(\mb_s,\beta)$ is {\em good enough} (with respect to $\alpha$),
if it is not flawless, but 
$\dom(\beta)=\dom(\alpha)$, and (\ref{ehto1}) holds 
for all $j\in\dom(\alpha)$ and all $r\not=s$, and there is a string
$t\in\{0,1\}^n$ with $|\{i:s_i\not=t_i\}|=1$ such that  
\begin{equation}\label{ehto2}
\alpha(j)=a_s\iff\beta(j)=c_t
\end{equation}
for all $j\in\dom(\alpha)$. Thus, if  $a_s\in\ran(\alpha)$,
then for each $t\in\{0,1\}^n$ with $|\{i:s_i\not=t_i\}|=1$
there is a unique $\beta$ such that $(\mb_s,\beta)$ is
good enough; we denote this $\beta$ by $\beta_{s,t,\alpha}$.

For the rest of this section, $\B$ will always denote a set
of structures of the form $(\mb_s,\beta)$, and 
$\A$ will denote a singleton set $\{(\ma,\alpha)\}$.
The {\em measure} of the set $\B$ is defined to be
$$
  M(\B)=(n+1)\cdot f(\B)+g(\B),
$$
where $f(\B)$ is the number of flawless structures
in $\B$ and $g(\B)$ is the number of good enough structures
in $\B$.

\begin{lemma}\label{atomi}
If $M(\B)>1$, then there is no atomic formula
$\p$ such that $(\A,\B)\models\p$ or $(\A,\B)\models\neg\p$.
\end{lemma}

\proof 
If $M(\B)>1$, then either there is a flawless structure in $\B$,
or there are at least two good enough structures in $\B$. 
If $(\mb_s,\beta)\in\B$
is flawless, then by condition (\ref{ehto1}), $(\ma,\alpha)$ and 
$(\mb_s,\beta)$
satisfy the same atomic formulas, whence no atomic formula separates
$\A$ and $\B$. 

Assume then, that $(\mb_s,\beta_{s,t,\alpha})$ and
$(\mb_{s'},\beta_{s',t',\alpha})$ are two distinct good enough structures 
in $\B$.
Let $\p$ be an atomic formula. If $\p$ is an identity $x_j=x_k$
(with $j,k\in\dom(\alpha)$), then it follows easily from (\ref{ehto1})
and (\ref{ehto2}) that $(\ma,\alpha)\models\p\iff(\mb_s,\beta_{s,t,\alpha})
\models\p$.
Thus, $\p$ does not separate $\A$ and $\B$.

Consider then the case $\p=P_l(x_j)$, where $j\in\dom(\alpha)$.
Let $\alpha(j)=a_r$. If $r\not=s$, then $\beta_{s,t,\alpha}(j)=b_r$, and we
have
$$(\ma,\alpha)\models\p\iff r_l=1\iff(\mb_s,\beta_{s,t,\alpha})\models\p.$$ 
Similarly,
if $r\not=s'$, then 
$(\ma,\alpha)\models\p\iff(\mb_{s'},\beta_{s',t',\alpha})\models\p$.
Assume finally, that $r=s=s'$. 
Since $(\mb_s,\beta_{s,t,\alpha})\not=(\mb_{s'},\beta_{s',t',\alpha})$, 
we have $t\not=t'$. Moreover, since $|\{i:r_i\not=t_i\}|=
|\{i:r_i\not=t'_i\}|=1$, either $r_l=t_l$ or $r_l=t'_l$. Thus, it is
not possible that $\p$ separates 
$(\ma,\alpha)$ from both $(\mb_s,\beta_{s,t,\alpha})$ and 
$(\mb_{s'},\beta_{s',t',\alpha})$.
We conclude that in all cases 
$(\A,\B)\not\models\p$ and $(\A,\B)\not\models\neg\p$.
\qed

\begin{lemma}\label{askel}
(a) If $\B=\C\cup\D$, then $M(\C)+M(\D)\ge M(\B)$.

(b) If $\A'=\A(F/j)$ and $\B'=\B(\star/j)$, then
$M(\B')\ge M(\B)-1$.
\end{lemma}

\proof (a) If $\B=\C\cup\D$, then obviously $f(\C)+f(\D)\ge f(\B)$, and
$g(\C)+g(\D)\ge g(\B)$. Hence we have  $M(\C)+M(\D)= (n+1)
(f(\C)+f(\D))+(g(\C)+g(\D))
\ge M(\B)$.

(b) Let $F((\ma,\alpha))=a_r$.
Thus, $\A'=\{(\ma,\alpha')\}$, where $\alpha'=\alpha(a_r/j)$.
Observe first that if $(\mb_s,\beta_{s,\alpha})$ is a flawless structure, 
and $r\not=s$, then $(\mb_s,\beta_{s,\alpha}(b_r/j))=
(\mb_s,\beta_{s,\alpha'})$ is also flawless, and clearly
$(\mb,\beta_{s,\alpha})\in\B\iff (\mb,\beta_{s,\alpha'})\in\B'$.

Assume then that $(\mb_s,\beta_{s,t,\alpha})$ is a good enough structure. 
If $r\not=s$, then $(\mb_s,\beta_{s,t,\alpha}(b_r/j))=
(\mb_s,\beta_{s,t,\alpha'})$ is also good enough. On the other hand,
if $r=s$, then $(\mb_s,\beta_{s,t,\alpha}(c_t/j))=
(\mb_s,\beta_{s,t,\alpha'})$ is good enough. In both cases,
$(\mb,\beta_{s,t,\alpha})\in\B\iff (\mb,\beta_{s,t,\alpha'})\in\B'$.

Thus, we see that if $\B$ does not contain a flawless structure 
of the form $(\mb_r,\beta_{r,\alpha})$, then $f(\B')=f(\B)$ and
$g(\B')=g(\B)$, whence the claim is true.

Assume finally, that there is a flawless structure 
$(\mb_r,\beta_{r,\alpha})$ in $\B$. Since $\alpha'(j)=a_r$, no
structure $(\mb_r,\beta)$ is flawless with respect to $\alpha'$.
On the other hand, for each $t\in\{0,1\}^n$ with $|\{i:r_i\not=t_i\}|=1$,
there is a new good enough structure $(\mb_r,\beta_{r,\alpha}(c_t/j))=
(\mb_r,\beta_{r,t,\alpha'})$ in $\B'$.
Thus, in this case $f(\B')=f(\B)-1$ and $g(\B')=g(\B)+n$, whence
$M(\B')= (n+1) (f(\B)-1)+g(\B)+n=M(\B)-1$.
\qed

\begin{lemma}\label{voitto}
If $w<M(\B)$, then player II has a winning strategy in 
$\ef^\exists_w(\A,\B)$.
\end{lemma}

\proof We prove the claim by induction on $w$. Consider first the case
$w=1$. By the definition of the game $\ef^\exists_1(\A,\B)$, there 
are no moves, and player I wins only if there is an atomic formula
$\p$ such that $(\A,\B)\models\p$ or $(\A,\B)\models\neg\p$. Since 
$M(\B)>w=1$, by Lemma~\ref{atomi}, there is no such $\p$.

Assume then that $w>1$,
and the claim is true for all $u<w$. Since $M(\B)>w\ge 1$, by 
Lemma~\ref{atomi} again, player I does not win the game without making
moves. Consider then the options of player I for his first move. 
Making a left splitting move $\A=\C\cup\D$ is not possible, since
$\A$ is a singleton $\{(\ma,\alpha)\}$. 
Suppose then that player I makes a right splitting move $w=u+v$
and $\B=\C\cup\D$. Then by Lemma~\ref{askel}(a), $M(\C)+M(\D)\ge
M(\B)$, and since $w<M(\B)$, either $u<M(\C)$, or $v<M(\D)$.
If $u<M(\C)$, then by induction hypothesis, player II has 
a winning strategy in the game $\ef^\exists_u(\A,\C)$. Similarly,
if $v<M(\D)$, then player II has 
a winning strategy in the game $\ef^\exists_v(\A,\D)$. Thus, by
choosing the appropriate position $(u,\A,\C)$ or $(v,\A,\D)$, player
II is guaranteed to win the game $\ef^\exists_w(\A,\B)$.

Suppose then that player I starts with a left supplementing move
$j$ and $F$, where $F$ is a choice function for $\A$. 
The next position in the game is then $(w-1,\A',\B')$, where
$\A'=\A(F/j)$ and $\B'=\B(\star/j))$. By Lemma~\ref{askel}(b)
and our assumption $w<M(\B)$, we have $w-1<M(\B)-1\le M(\B')$,
whence by induction hypothesis, player II has a winning
startegy in the continuation of the game $\ef^\exists_w(\A,\B)$
from position $(w-1,\A',\B')$ onwards.
\qed

Consider now the classes $\A_0$ and $\B_0$ defined in the beginning
of this section. Since the variable assignment in the only structure
in $\A_0$ is empty, all the $2^n$ structures $(\mb_s,\emptyset)$
in $\B_0$ are flawless. Thus, $M(\B_0)= (n+1)f(\B_0)+g(\B_0)=(n+1)2^n$, and 
by Lemma~\ref{voitto}, player II has a winning strategy in the game
$\ef^\exists_w(\A_0,\B_0)$ whenever $w<(n+1)2^n$. As all Boolean combinations
of the predicates $P_1,\ldots,P_n$ are non-empty in $\ma$, but each
structure in $\B_0$ has an empty Boolean combination, we get the desired
lower bound result:

\begin{corollary}
If $\p$ is an existential first order sentence expressing 
the property that all Boolean combinations of $n$ unary predicates
are non-empty, then the size of $\p$ is at least $(n+1)2^n$. \qed
\end{corollary}

\section{The Existential Complexity
of the Length of Linear Order}\label{order}

As we noted in the introdution, for each $n$ there is first
order sentence $\p_n$ of logarithmic size with respect to $n$
which expresses the property that
the length of a linear order is at least $n$. However, the sentence $\p_n$
has an unbounded number of quantifier alternations. In this section we
show that $2n-1$ is the minimum size of an existential sentence
expressing this property. 

To prove the upper bound, define the following sequence of existential
formulas:
\begin{itemize}
\item[] $\s_2:= \exists x_1\exists x_2 (x_1<x_2)$, \quad and
\item[] $\s_{n+1}:= \exists x_{n+1}(\s_n\land x_n<x_{n+1})$ \quad 
for all $n\ge 2$.
\end{itemize}
Clearly $\s_n$ is true in a linear order if and only if its length
is at least $n$, and an easy induction shows that $\rk(\s_n)=2n-1$.

To prove the lower bound, we will use again the existential game 
$\ef^\exists_w$.
Let $\A_0=\{(\ma,\emptyset)\}$, where $\ma$ is a linear order of length $n$, 
and let $\B_0=\{(\mb,\emptyset)\}$, where $\B$ is a linear order of 
length $n-1$. Our aim is to show that player II has a winning strategy
in the game $\ef^\exists_w(\A_0,\B_0)$ for all $w<2n-1$. 

Consider a position $(u,\A,\B)$ in the game $\ef^\exists_w(\A_0,\B_0)$.
Since the game is existential, $\A$ consists of a single structure
$(\ma,\alpha)$. Let $a_1<^\ma\cdots<^\ma a_l$ be the elements in 
$\ran(\alpha)$, and let $a_0$ and $a_{l+1}$ be the least and the largest
element in $\ma$, respectively.
We say that a variable assignment $\beta$ in $\mb$ is {\em acceptable} (with
respect to $\alpha$) if $\dom(\beta)=\dom(\alpha)$, there are elements
$b_1\le^\mb\cdots\le^\mb b_l$ such that $\ran(\beta)=\{b_1,\ldots,b_l\}$
and for all $i\in\{1,\ldots,l\}$ and all $j\in\dom(\alpha)$
\begin{equation}\label{ehto3}
\alpha(j)=a_i\iff\beta(j)=b_i.
\end{equation}
In other words, $\beta$ is acceptable if and only if
the mapping $\alpha(j)\mapsto\beta(j)$, $j\in\dom(\alpha)$,
preserves the relation $\le$.

Furthermore, we say that $\beta$ is {\em nice} 
(with respect to $\alpha$),
if it is acceptable, and in addition $|\{i\le l:d(a_i,a_{i+1})
\not=d(b_i,b_{i+1})\}|=1$, where $d(x,y)$ is the distance
between $x$ and $y$ in the given linear order, and 
$b_0$ and $b_{l+1}$ are 
the least and the largest element in $\mb$, respectively. The 
{\em distance from defect} of $\beta$ is defined 
by $\delta(\beta)=d(b_i,b_{i+1})$, where $i\le l$ is the unique
index such that $d(a_i,a_{i+1})\not=d(b_i,b_{i+1})$;
we denote this index by $i(\beta)$.
Note that $d(b_i,b_{i+1})=d(a_i,a_{i+1})-1$ for $i=i(\beta)$.
Note also, that for each $i\le l$ there is exactly one 
nice $\beta$ such that $i(\beta)=i$.

For the rest of this section, we assume that $\A$ is a singleton
set $\{(\ma,\alpha)\}$ and $\B$ is a set of structures
of the form $(\mb,\beta)$. The {\em niceness measure} of
$\B$ is defined to be
$$
  N(\B)=\sum_{\beta\in\N} (2\delta(\beta)+1),
$$
where $\N$ is the set of all nice variable assignments
$\beta$ such that $(\mb,\beta)\in\B$.

\begin{lemma}\label{perus}
If $N(\B)>1$, then there is no atomic formula
$\p$ such that $(\A,\B)\models\p$ or $(\A,\B)\models\neg\p$.
\end{lemma}

\proof If $N(\B)>1$, then either there is a nice assignment $\beta\in\N$
such that $\delta(\beta)\ge 1$, or there are two distinct nice
assignments $\beta,\beta'\in\N$. Assume first that $\beta$ is a nice
assignment in $\N$, and $\delta(\beta)\ge 1$. Then there are
elements $a_1,\ldots,a_l,b_1,\ldots,b_l$ such that $\ran(\alpha)=
\{a_1,\ldots,a_l\}$, $\ran(\beta)=\{b_1,\ldots,b_l\}$, 
$a_1<^\ma\cdots<^\ma a_l$ and $b_1\le^\mb\cdots\le^\mb b_l$. 
Since $\beta$ is nice, $d(b_i,b_{i+1})=\delta(\beta)\ge 1$ for
$i=i(\beta)$, and $d(b_i,b_{i+1})=d(a_i,a_{i+1})\ge 1$
for all other $i\le l$. Thus, we have $b_1<^\mb\cdots<^\mb b_l$.
It follows easily from condition (\ref{ehto3}) that $(\ma,\alpha)$
and $(\mb,\beta)$ satisfy the same atomic formulas, whence
no atomic formula separates $\A$ and $\B$.

Assume then that $\beta,\beta'\in\N$, and $\beta\not=\beta'$.
As noted above, this means that $i(\beta)\not=i(\beta')$.
Let $a_1,\ldots,a_l$ be as above, and let $\ran(\beta)=
\{b_1,\ldots,b_l\}$ and $\ran(\beta')=\{b'_1,\ldots,b'_l\}$ with
$b_1\le^\mb\cdots\le^\mb b_l$ and $b'_1\le^\mb\cdots\le^\mb b'_l$.
As above, we see that $b_i<^\mb b_{i+1}$ for all $i\le l$ except
$i=i(\beta)$, and similarly $b'_i<^\mb b'_{i+1}$ for all $i\le l$ except
$i=i(\beta')$. Moreover, since $i(\beta)\not=i(\beta')$, for all
$i,j\le l$, we have either ($a_i<^\ma a_j\iff b_i<^\mb b_j$) or
($a_i<^\ma a_j\iff b'_i<^\mb b'_j$). Using condition (\ref{ehto3}),
it is now easy to see that no atomic formula separates the sets
$\A$ and $\B$.
\qed

\begin{lemma}\label{seur}
(a) If $\B=\C\cup\D$, then $N(\C)+N(\D)\ge N(\B)$.

(b) If $\A'=\A(F/j)$ and $\B'=\B(\star/j)$, then
$N(\B')\ge N(\B)-1$.
\end{lemma}

\proof (a) Assume that $\B=\C\cup\D$. Let $\N$, $\O$ and $\P$ be
the sets of nice assignments $\beta$ such that $(\mb,\beta)\in\B$,
$(\mb,\beta)\in\C$ and $(\mb,\beta)\in\D$, respectively.
Then we have

\begin{eqnarray*}
N(\B)&=&\sum_{\beta\in\N} (2\delta(\beta)+1) \\ 
&\le&\sum_{\beta\in\O} (2\delta(\beta)+1)+\sum_{\beta\in\P} 
(2\delta(\beta)+1)=N(\C)+N(\D).
\end{eqnarray*}

(b) Let $\N$ and $\N'$ be the sets of nice $\beta$
such that $(\mb,\beta)$ is in $\B$ and $\B'$, respectively.
Let $\ran(\alpha)=\{a_1,\ldots,a_l\}$ with
$a_1<^\ma\cdots<^\ma a_l$, where $\alpha$ is the assignment
such that $\A=\{(\ma,\alpha)\}$. As usual, we denote by
$a_0$ and $a_{l+1}$ the least and the largest element in $\ma$, 
respectively. Furthermore, we denote $F((\ma,\alpha))$ by $c$.

Assume first that $c=\alpha(k)$ for some $k\in\dom(\alpha)$. 
It is easy to see that for any assignment $\beta$, 
$$
\beta\in\N\iff\beta(d/j)\in\N',
$$ 
where $d=\beta(k)$, and moreover $\delta(\beta(d/j))=\delta(\beta)$. 
Note also that if $\beta\not=\beta'$,
then $\beta(d/j)\not=\beta'(d'/j)$. Thus, in this case we have
$N(\B')=N(\B)$.

Assume next that $c=a_0$ or $c=a_{l+1}$. Then as above, we see
that $\beta\in\N\iff\beta(d/j)\in\N'$, where $d$ is the least or the
largest element in $\B$, respectively, and 
$\delta(\beta(d/j))=\delta(\beta)$. Thus, also in this case we
conclude that $N(\B')=N(\B)$.

Assume finally, that $c\not\in\{a_0,\ldots,a_{l+1}\}$. Then there 
is an index $h\le l$ such that $a_h<^\ma c<^\ma a_{h+1}$. 
Let $\beta$
be a nice assignment in $\N$, and let $\ran(\beta)=\{b_1,\ldots,b_l\}$
with $b_1\le^\mb\cdots\le^\mb b_l$. 
As noted earlier,
there is exactly one nice assignment $\beta_h$ such that $i(\beta_h)=h$. 
If $\beta\not=\beta_h$, then 
$d(b_h,b_{h+1})=d(a_h,a_{h+1})$, whence there is an element $d$
such that $d(b_h,d)=d(a_h,c)$ and $d(d,b_{h+1})=d(c,a_{h+1})$.
This means that $\delta(\beta(d/j))=\delta(\beta)$, and clearly
$\beta(d/j)$ is in $\N'$. 

On the other hand, if $\beta=\beta_h$, then 
$d(b_h,b_{h+1})=d(a_h,a_{h+1})-1$, and there are elements
$d$ and $e$ such that $d(b_h,d)=d(a_h,c)-1$, $d(d,b_{h+1})=d(c,a_{h+1})$,
$d(b_h,e)=d(a_h,c)$ and $d(e,b_{h+1})=d(c,a_{h+1})-1$. Let
$\beta'=\beta_h(d/j)$ and $\beta''=\beta_h(e/j)$. Then $\beta',\beta''\in\N'$,
and we have
$$
  \delta(\beta')+\delta(\beta'')=d(b_h,d)+d(e,b_{h+1})=
  d(a_h,a_{h+1})-2=\delta(\beta_h)-1,
$$
whence
$$
  (2\delta(\beta')+1)+(2\delta(\beta'')+1)=2(\delta(\beta_h)-1)+2=
  (2\delta(\beta_h)+1)-1.
$$
Thus, if $\beta_h$ is in $\N$, we get
\begin{eqnarray*}
  N(\B')&=&\sum_{\beta'\in\N'}(2\delta(\beta')+1)\\ &=&
  \sum_{\beta\in\N\setminus\{\beta_h\}}(2\delta(\beta)+1)
  +((2\delta(\beta_h)+1)-1)=N(\B)-1,
\end{eqnarray*}
and if $\beta_h$ is not in $\N$, we have $N(\B')=N(\B)$.
\qed

\begin{lemma}\label{vstrat}
If $w<N(\B)$, then player II has a winning strategy in 
$\ef^\exists_w(\A,\B)$.
\end{lemma}

\proof The proof is {\em verbatim} the same as for Lemma~\ref{voitto};
just replace $M(\B)$ by $N(\B)$.
\qed

The proof of the lower bound result is also similar to that in the
previous section. Consider the sets $\A_0=\{(\ma,\emptyset)\}$
and $\B_0=\{(\mb,\emptyset)\}$. Clearly $\delta(\emptyset)=n-1$,
and so $N(\B_0)=2(n-1)+1=2n-1$. Thus, by Lemma~\ref{vstrat}, player II
has a winning strategy in the game $\ef^\exists_w(\A_0,\B_0)$ for
all $w<2n-1$. Since the length of the linear order $\ma$ is $n$,
while the length of $\mb$ is less than $n$, we have proved

\begin{corollary}
If $\p$ is an existential first order sentence expressing 
the property that the length of a linear order is at least $n$, 
then the size of $\p$ is at least $2n-1$. \qed
\end{corollary}

\bibliographystyle{plain}

\end{document}